\documentclass{amsart}
\usepackage[utf8]{inputenc}
\usepackage{pifont}
\usepackage[numbers]{natbib}
\usepackage{enumerate,verbatim,calc}
\usepackage{braket}
\usepackage{graphicx}
\usepackage[top=1.2in, bottom=1.2in, left=1.4in, right=1.4in]{geometry}
\usepackage[mathscr]{euscript}
\usepackage{mathrsfs}
\usepackage{amsmath,euscript,amsfonts,amsthm,amssymb,latexsym}

\def\R{{\mathbb R}}
\def\A{{\mathbb A}}
\def\C{{\mathbb C}}

\def\B{{\mathbb B}}

\def\T{{\mathbb T}}

\def\O{{\mathbb O}}
\def\A{{\mathbb A}}

\def\h2{{\mathcal H}^2(\B^{2m})}
\def\b2{{\mathcal A}^2(\B^{2m})}

\def\L{{\mathcal L}}
\def\T{{\mathcal T}}

\begin{document}
\title{\bf Multivariable Isometries Related to Certain Convex Domains}

\author{Ameer Athavale}
\address{Department of Mathematics, Indian Institute of Technology Bombay, Powai, Mumbai 400076, India}
\email{athavale@math.iitb.ac.in}
\begin{abstract}
There exist several interesting results in the literature on subnormal operator tuples having their spectral properties tied to the geometry of strictly pseudoconvex domains or to that of bounded symmetric domains in $\C^n$. We introduce a class $\Omega^{(n)}$ of convex domains in $\C^n$ which, for $n \geq 2$, is distinct from the class of strictly pseudoconvex domains and the class of bounded symmetric domains and which lends itself for the application of the theories related to the abstract inner function problem and the $\bar \partial$-Neumann problem, allowing us to make a number of interesting observations about certain subnormal operator tuples associated with the members of the class $\Omega^{(n)}$.
\end{abstract}
\maketitle
\section{Introduction}
We use ${\mathcal B}({\mathcal H})$ to denote the algebra of bounded linear operators on  a complex infinite-dimensional separable Hilbert space $\mathcal H$ and use $I$ to denote the identity operator on $\mathcal H$. 
An $n$-tuple $S = (S_1,\ldots, S_n)$ of commuting operators $S_i$ in ${\mathcal B}({\mathcal H})$ is said to be {\it subnormal} if there exist a Hilbert
space ${\mathcal K}$ containing ${\mathcal H}$ and an $n$-tuple $N = (N_1, \ldots,
N_n)$ of commuting normal operators $N_i$ in ${\mathcal B}({\mathcal K})$ such
that $N_i {\mathcal H}  \subset {\mathcal H}$ and $N_i\vert{\mathcal H} = S_i$ for $1
\leq i \leq n$.  Among all the normal extensions of a subnormal tuple $S$, there is a `minimal normal extension' which is unique up to unitary equivalence (see \cite{I}). An $n$-tuple $T = (T_1,\ldots, T_n)$ of commuting operators $T_i$ in ${\mathcal B}({\mathcal H})$ is said to be {\it essentially normal} if the operators $T_i^*T_j -T_jT_i^*$ are compact for all $i$ and $j$, while $T$ is said to be {\it cyclic} if there exists a vector $f$ in ${\mathcal H}$ (referred to as a {\it cyclic vector} for $T$) such that the linear span $\vee\{T_1^{k_1}T_2^{k_2}\cdots T_n^{k_n}f: k_i$ are non-negative  integers$\}$ is dense in ${\mathcal H}$. There exist several interesting results in the literature on subnormal operator tuples (and in particular on  essentially normal and/or cyclic subnormal operator tuples) having their spectral properties tied to the geometry of strictly pseudoconvex domains or to that of bounded symmetric domains in $\C^n$ (refer, for example, to \cite{At3}, \cite{At-Po}, \cite{Di1}, \cite{Di2}, \cite{Di-E}, \cite{Di-E-Ev}, \cite{E}, \cite{E-Ev}, \cite{U}). These results are largely manifestations of the functional calculus for subnormal operator tuples thriving upon some elegant function-theoretic results valid in the context of those two types of domains. (We refrain from referrring to an endless list of papers that specifically deal with subnormal operator tuples related to the unit ball ${\mathbb B}_n$ in $\C^n$ which is a strictly pseudoconvex as well as a bounded symmetric domain).\\

In Section 2 we introduce a class $\Omega^{(n)}$ of convex domains in $\C^n$ whose members  are parametrized by $n$-tuples $p$ with the coordinates of $p$ being tuples (of varying lengths) of positive integers subject to certain constraints. For $n \geq 2$, the class $\Omega^{(n)}$ of domains $\Omega_p$ turns out to be distinct from the class of strictly pseudoconvex domains and the class of bounded symmetric domains. The new class allows for the application of the theory related to the abstract inner function problem (refer to \cite{A} and \cite{Di-E}) as well as of the theory related to the  $\bar \partial$-Neumann problem (refer to \cite{Bo-S} and \cite{F-Ko}). The multiplication tuples associated with the Hardy-type function spaces associated with the domains $\Omega_p$ turn out to be so-called (regular) $A$-isometries. We record a few properties of the domains $\Omega_p$ that are relevant for the application of some known results in the literature to those $A$-isometries; these applications mostly result from the existence of an abundance of inner functions on the domains $\Omega_p$ as in the case of domains that are either strictly pseudoconvex or bounded symmetric (refer to \cite{A} and \cite{Di-E}).\\

In Section 3 we record parts of the theory related to the $\bar \partial$-Neumann problem and the tangential Neumann problem as are of interest to us. The $\bar \partial$-Neumann problem (resp. tangential  Neumann problem) will be seen to be of particular relevance in the context of the multiplication tuples $M_{\nu_p,z}$ (resp. $M_{\sigma_p,z})$ associated with the Bergman (resp. Hardy) spaces of the domains $\Omega_p$. Indeed, among our  concerns in Section 3 will be the compactness of the so-called $\bar \partial$-Neumann operator and that of the so-called tangential Neumann operator, since the compactness of the $\bar \partial$-Neumann operator (resp. tangential Neumann operator) guarantees the essential normality of the tuple $M_{\nu_p,z}$ (resp. $M_{\sigma_p,z})$. \\

In Section 4 we discuss multivariable isometries associated with certain convex domains $\Sigma_p$ that are more general than the domains $\Omega_p$,  providing an intrinsic characterization of such multivariable isometries (referred to as $\partial \Sigma_p$-isometries). In particular, a succinct characterization of a  $\partial \Sigma_p$-isometry is derived for a special type of $\Sigma_p$, which is an apt generalization of that of a `spherical isometry' (see \cite{At2}). We also dwell there on the intertwining of a $\partial \Omega_p$-isometry with certain other subnormal tuples. Finally, we elaborate upon the significance of the domains $\Omega_p$ for some operator theoretic considerations that go beyond the topic of multivariable isometries. \\

For any terminology employed from the area of several complex variables and for any standard results quoted from there, the references \cite{J-P}, \cite{Kr} and \cite{Ran} should be more than adequate.  
\section{Convex domains $\Omega_p$}
Let $p= (p_1,p_2,\ldots,p_n)$ be an $n$-tuple of $m_i$-tuples $p_i= (p_{i,1},p_{i,2},\ldots,p_{i,m_i})$ where, for each $i$ satisfying $1 \leq i \leq n$, $ p_{i,1}$, $p_{i,2}$,...,$p_{i,m_i}$ (with $m_i \geq 2$) are relatively prime positive integers so that ${\rm gcd}\{p_{i,1},p_{i,2},\ldots,p_{i,m_i}\}=1$. The subset $\Omega_p$ of $\C^n$ is defined by $\Omega_p = \{z =(z_1,z_2,\ldots,z_n) \in \C^n: 
\sum_{i=1}^n\sum_{j=1}^{m_i} |z_i|^{2p_{i,j}} < 1\}$. The set $\Omega_p$ is easily seen to be a convex complete Reinhardt domain in $\C^n$ with the real analytic boundary $\partial\Omega_p=\{z =(z_1,z_2,\ldots,z_n) \in \C^n: \sum_{i=1}^n\sum_{j=1}^{m_i} |z_i|^{2p_{i,j}} =1\}$. Some of the results in Sections 2 and 3 as stated for the domains $\Omega_p$ also hold for certain domains more general than $\Omega_p$ - these will be pointed out explicitly in Section 4. We use the symbol $\Omega^{(n)}$ to denote the class of domains $\Omega_p$ in $\C^n$ parametrized by the tuples $p$ as described above. For $z \in \C$, $\bar z$ denotes the complex conjugate of $z$ and, for any complex-valued function $\phi$, $\bar \phi$ is the function satisfying $\bar \phi(z) = \overline{\phi(z)}$. \\

{\bf Remark 2.1}. For $n=1$, the domains $\Omega_p$ reduce to the open unit disks in the plane (of various radii) centered at the origin for which the theme of the paper stands already well-explored (refer, for example, to \cite{Co1} and \cite{D}). For that reason, and for the validity of certain assertions to follow, {\bf we assume hereafter in any discussion involving $\Omega_p$ that $n \geq 2$}. \\

{\bf Remark 2.2}. The domain $\Omega_p$ equals $\{z \in \C^n: u(z) < 0\}$ where  $u(z) =  \sum_{i=1}^n\sum_{j=1}^{m_i} |z_i|^{2p_{i,j}} - 1$. For $b \in \partial \Omega_p$, let $\T_b(\partial \Omega_p) = \{X=(X_1,\ldots,X_n)\in \C^n: \sum_{j=1}^n \frac{\partial u}{\partial z_j}(b)X_j =0\}$ be the complex tangent space to $\partial\Omega_p$ at $b$. The Levi form $\L u(b,X) = \sum_{j,k =1}^n\frac{\partial^2u}{\partial z_j \partial {\bar z}_k}(b) X_j {\bar X}_k$ is non-negative for every $b \in \partial\Omega_p$ and $X \in \T_b(\partial \Omega_p)$. However, for not all permissible choices of $p$, the Levi form $\L u(b,X)$ is positive for every $b \in \partial\Omega_p$ and every non-zero $X \in \T_b(\partial \Omega_p)$. Thus $\Omega_p$ (though a pseudoconvex domain) is not in general strictly pseudoconvex at every point of its boundary $\partial \Omega_p$, and the class $\Omega^{(n)}$ is distinct from the class of strictly pseudoconvex domains in $\C^n$. \\ 

{\bf Remark 2.3}. By a result of Cartan \cite{Car}, every bounded symmetric domain $D$ in $\C^n$ is homogeneous in the sense that the automorphism group of $D$ acts transitively on $D$. Also, by a result of Pinchuk \cite{Pi}, every bounded homogeneous domain in $\C^n$ with smooth boundary is biholomorphically equivalent to the unit ball ${\mathbb B}_n$ in $\C^n$. If $\Omega_p$ were to be a bounded symmetric domain, it would thus be  biholomorphically equivalent to ${\mathbb B}_n$. A result of Sunada \cite{Su}, however, states that two Reinhardt domains $D_1$ and $D_2$ in $\C^n$ that contain the origin are biholomorphically equivalent if and only if there exist positive numbers $r_1,\ldots,r_n$ and a permutation $\sigma$ of $\{1,\ldots,n\}$ such that $D_2 = \{(r_1z_{\sigma(1)},\ldots,r_nz_{\sigma(n)}): (z_1,\ldots,z_n) \in D_1 \}$. It follows that the class $\Omega^{(n)}$ is distinct from the class of bounded symmetric domains in $\C^n$. \\ 

Let $K \subset \C^n$ be compact, and let $A$ be a unital closed subalgebra of $C(K)$ containing $n$-variable complex polynomials. The {\it Shilov boundary} of $A$ is defined to be the smallest closed subset $S$ of $K$ such that
$$ \|f\|_{\infty,K} = \|f\|_{\infty,S}\ \ \ (f \in A).$$
Of special interest to us is the subalgebra $A(\Omega_p) = \{f \in C({\bar {\Omega}_p}): f {\rm \ is \ holomorphic \ on\ } \Omega_p\}$ of $C({\overline \Omega}_p)$, where  ${\overline {\Omega}_p} = \{z =(z_1,z_2,\ldots,z_n) \in \C^n: \sum_{i=1}^n\sum_{j=1}^{m_i} |z_i|^{2p_{i,j}} \leq 1\}$ is the closure of $\Omega_p$. If $\O({\overline {\Omega}_p})$ is the vector space of functions $f$ such that $f$ is holomorphic on an open neighborhood $U_f$ of $\bar \Omega_p$, then (referring to the first line of Remark 3.2) it is easy to see that the closure of $\O({\overline {\Omega}_p})$ in the sup norm with respect to $\bar \Omega_p$ is $A(\Omega_p)$.\\

{\bf Proposition 2.4}. The Shilov boundary of $A(\Omega_p)$ coincides with the topological boundary $\partial\Omega_p$ of $\Omega_p$.
\begin{proof} Since $\Omega_p$ is a bounded pseudoconvex domain in $\C^n$ with smooth boundary, it follows from \cite[Folgerung 5]{Pf} (see also \cite{H-S}) that the Shilov boundary of $A(\Omega_p)$ is the closure of the set of  strictly pseudoconvex points in $\partial \Omega_p$.  It is easy to see that any point $b=(b_1,\cdots,b_n)$ of $\partial \Omega_p$ for which each $b_i$ is non-zero is a point of strict pseudoconvexity. But such points are dense in $\partial\Omega_p$ so that the Shilov boundary of  $A(\Omega_p)$ is $\partial\Omega_p$. 
\end{proof}

Let $K$ be a compact subset of $\C^n$, let $A$ be a closed subspace of $C(K)$, and let $\eta$ be a positive regular Borel measure on $K$. The triple $(A,K,\eta)$ is said to be {\it regular} (in the sense of \cite{A}) if, for any positive function $\phi$ in $C(K)$, there exists a sequence of functions $\{\phi_m\}_{m \geq 1}$ in $A$ such that $|\phi_m| < \phi$ on $K$ and $\lim_{m \rightarrow \infty} |\phi_m| = \phi$ $\eta$-almost everywhere. \\

{\bf Proposition 2.5}. For any positive regular Borel measure $\mu_p$ on ${\overline {\Omega}_p}$ with supp$(\mu_p) \subset \partial\Omega_p$, the triple $(A(\Omega_p), {\overline {\Omega}_p}, \mu_p)$ is regular as is the triple $(A(\Omega_p)|\partial\Omega_p, \partial\Omega_p, \mu_p)$.
\begin{proof}
For $p= (p_1,p_2,\ldots,p_n)$ with $p_i= (p_{i,1},p_{i,2},\ldots,p_{i,m_i})$, let $N = m_1+\cdots+m_n$. 
Consider $f(z) = ({z_1}^{p_{1,1}},\ldots,{z_1}^{p_{1,m_1}},\ldots,{z_n}^{p_{n,1}},\ldots,{z_n}^{p_{n,m_n}})$, $z \in {\overline {\Omega}_p}$. Clearly, $f$ maps $\partial\Omega_p$ into the topological
 boundary $\partial\mathbb B_{N}$ of the unit ball $\mathbb B_{N}$ of $\C^N$. Thus the regularity of the triple $(A(\Omega_p), {\overline {\Omega}_p}, \mu)$ will follow from \cite [Proposition 2.5]{Di-E} provided we verify $f$ to be injective. The regularity of the triple $(A(\Omega_p)|\partial\Omega_p, \partial\Omega_p, \mu)$ will then be an easy consequence of the Tietze extension theorem. Thus let $z=(z_1,\ldots,z_n)$ and $w=(w_1,\ldots,w_n)$ be distinct points of ${\overline {\Omega}_p}$ so that $z_i \neq w_i$ for some $i$. If only one of $z_i$ and $w_i$ is non-zero, then clearly
 $f(z) \neq f(w)$. So suppose both $z_i$ and $w_i$ are non-zero. Using the coprimality of $p_{i,1},p_{i,2},\ldots,p_{i,m_i}$, we choose integers $n_{i,1},n_{i,2},\ldots,n_{i,m_i}$ such that $\sum_{j=1}^{m_i}n_{i,j}p_{i,j} =1$. If one were to have ${z_i}^{p_{i,j}} = {w_i}^{p_{i,j}}$ for every $j$ such that $1 \leq j \leq m_i$, then that would clearly force the contradiction $z_i = w_i$. Thus ${z_i}^{p_{i,j}} \neq {w_i}^{p_{i,j}}$ for some $j$ satisfying $1 \leq j \leq m_i$, showing that $f(z) \neq f(w)$.
\end{proof} 

At this stage we refer the reader to Section 2 of \cite{Di-E-Ev}. If $\mu$ is a scalar spectral measure of the minimal normal extension $N \in {\mathcal B}({\mathcal K})^n$ of a subnormal tuple $S \in {\mathcal B}({\mathcal H})^n$, then there is an isomorphism $\Psi_N$ of the von Neumann algebra $L^{\infty}(\mu)$ onto the von Neumann algebra $W^*(N)$ generated by $N_i \in {\mathcal B}({\mathcal K})$. The {\it restriction algebra} ${\mathcal R}_S = \{f \in L^{\infty}(\mu): \Psi_N(f){\mathcal H}\subset {\mathcal H}\}$ is a weak$^*$ closed subalgebra of $L^{\infty}(\mu)$.  Let $K \subset \C^n$ be compact and let $A$ be a unital closed subalgebra of $C(K)$ containing $n$-variable complex polynomials. Following \cite{E}, we call a subnormal tuple $S$  an {\it $A$-isometry} if the spectral measure of the minimal normal extension $N$ of $S$ is supported on the Shilov boundary  of $A$ and if $A$ is contained in  ${\mathcal R}_S$. Given a normalized positive regular Borel measure $\mu_p$ supported on $\partial \Omega_p$, we let $H^2(\mu_p)$ be the closure of $A(\Omega_p)$ in $\L^2(\mu_p)$. Letting $\sigma_p$ denote the normalized surface area measure on $\partial \Omega_p$, we refer to $H^2(\sigma_p)$ as the {\it Hardy space} of $\Omega_p$. In view of Proposition 2.4 and  in view of the discussion in Section 2 of \cite{Di-E-Ev}, the multiplication tuple $M_{\mu_p,z} = (M_{\mu_p,z_1},\ldots,M_{\mu_p,z_n})$ of multiplications by the coordinate functions $z_i$ on $H^2(\mu_p)$ is an $A(\Omega_p)$-isometry (and has the multiplication tuple $N_{\mu_p,z} = (N_{\mu_p,z_1},\ldots,N_{\mu_p,z_n})$ associated with $L^2(\mu_p)$ as its minimal normal extension); also, in the light of Proposition 2.5, $M_{\mu_p,z}$ is  {\it regular} in the sense of \cite{E} (that is, in the sense of  \cite[Definition 2.6]{Di-E-Ev}).\\

The preceding observations allow us to bring all the results in \cite{Di2}, \cite{Di-E-Ev}, \cite{E} and \cite{E-Ev} related to a regular $A$-isometry to bear upon the multiplication tuple $M_{\mu_p,z}$; we highlight in Remarks 2.6 and 2.7 below a few implications of the results in those references. We also point out that some of those results are derived exploiting Prunaru's work in \cite{Pr}.\\

{\bf Remark 2.6}. Let $P_{\mu_p}$ be the orthogonal projection of $\L^2(\mu_p)$ onto $H^2(\mu_p)$ and let, for $\phi \in L^{\infty}(\mu_p)$, $N_{\mu_p,\phi}$ denote the operator of multiplication by $\phi$ on $L^2(\mu_p)$. We let $T_{\mu_p,\phi}$ stand for $P_{\mu_p}N_{\mu_p,\phi}|H^2(\mu_p)$ and refer to $\T\left(M_{\mu_p,z}\right)=\{T_{\mu_p,\phi}: \phi \in L^{\infty}(\mu_p)\}$ as the set of {\it $M_{\mu_p,z}$-Toeplitz operators}. Also, we use $H^{\infty}_{A(\Omega_p)}(\mu_p)$ to denote the weak$^*$ closure  of $A(\Omega_p)$ in $L^{\infty}(\mu_p)$ and refer to any member $\theta$ of $H^{\infty}_{A(\Omega_p)}(\mu_p)$ satisfying $|\theta| =1$ $\mu_p$-almost everywhere as a {\it $\mu_p$-inner function}. It follows from  \cite[Corollary 3.3]{Di-E-Ev} that $\T\left(M_{\mu_p,z}\right)$ equals the set $\{X \in {\mathcal B}(H^2(\mu_p)): T_{\mu_p,\bar \theta}XT_{\mu_p,\theta} = X {\rm \ for \ every \ \mu_p-inner\ function\ \theta}\}$. Further, if $C^*\left(\T\left(M_{\mu_p,z}\right)\right)$ is the $C^*$-subalgebra of ${\mathcal B}(H^2(\mu_p))$ generated by $\T\left(M_{\mu_p,z}\right)$, ${\mathcal S}{\mathcal C}\left(M_{\mu,z}\right)$ the two-sided closed ideal in $C^*\left(\T\left(M_{\mu_p,z}\right)\right)$ generated by semicommutators $T_{\mu_p,\phi}T_{\mu_p,\psi}-T_{\mu_p,\phi\psi} \ (\phi, \psi \in L^{\infty}(\mu_p))$, and $\left(N_{\mu_p,z}\right)'$ the commutant in ${\mathcal B}(\L^2(\mu_p))$ of $\{N_{\mu_p,z_1},\ldots,N_{\mu_p,z_n}\}$, then Corollary 3.7 of \cite{Di-E-Ev} yields the existence of a short exact sequence of C$^*$-algebras
$$ 0 \rightarrow {\mathcal S}{\mathcal C}\left(M_{\mu_p,z}\right) \overset {\iota}\rightarrow C^*\left(\T\left(M_{\mu_p,z}\right)\right) \overset{\pi}\rightarrow \left(N_{\mu_p,z}\right)' \rightarrow 0$$
where $\iota$ is the inclusion map and $\pi$ is a unital $*$-homomorphism that is in fact a left inverse of the compression map $\rho: \left(N_{\mu_p,z}\right)' \rightarrow {\mathcal B}(H^2(\mu_p))$ given by $\rho(Y) = P_{\mu_p}Y|H^2(\mu_p)$, $Y \in \left(N_{\mu_p,z}\right)'$. \\

{\bf Remark 2.7}. (a) Let ${\mathcal A}_{M_{\mu_p,z}}$ be the weak$^*$-closed subalgebra of ${\mathcal B}(H^2(\mu_p))$ generated by $M_{\mu_p,z_i}$\ $(1 \leq i \leq n)$ and the identity operator on $H^2(\mu_p)$. It is a consequence of \cite[Corollary 6]{E} that the weak operator topology and the weak$^*$ operator topology coincide on ${\mathcal A}_{M_{\mu_p,z}}$ and that every unital weak$^*$-closed subalgebra of ${\mathcal A}_{M_{\mu_p,z}}$ is reflexive; in particular, $M_{\mu_p,z}$ is reflexive (refer to \cite{E} for the relevant definitions).\\
(b) It is a consequence of \cite[Corollary 2]{Di2} that the set $\T\left(M_{\mu_p,z}\right)$ of $M_{\mu_p,z}$-Toeplitz operators is $2$-hyperreflexive with the $2$-hyperreflexivity constant $\kappa_2\left(\T\left(M_{\mu_p,z}\right)\right)$ being less than or equal to $2$ (refer to \cite{K-Pt} for the relevant definitions).\\
 
As the results in \cite{Di-E-Ev} and \cite{E-Ev} show, one gets some extra mileage out of the notion of a regular $A$-isometry $T$ under the additional assumption that $T$ is essentially normal. We plan to explore the essential normality of the multiplication tuple $M_{\sigma_p, z}$ associated with the Hardy space $H^2(\sigma_p)$ of $\Omega_p$, and for that purpose we invoke in the next section the theory related to the famous $\bar \partial$-Neumann problem. 

\section{$\bar \partial$-Neumann Operator and the Tangential Neumann Operator}
While a basic reference for the material in this section is \cite{F-Ko}, we find in addition \cite{U} as a convenient reference for our purposes (see also \cite{S-Sh-U}). Indeed, some of the arguments in \cite{U} are adaptations and extensions of the arguments in \cite{F-Ko} to the context of the Hardy and Bergman spaces of strictly pseudoconvex domains and our task here is to push through the analogs of those arguments in the context of the domains $\Omega_p$. \\

Let $\Omega$ be a bounded pseudoconvex domain in $\C^n$ with its boundary $\partial \Omega = \{z \in \C^n : \rho(z)=0\}$ defined by a smooth function  $\rho : \C^n \rightarrow \R$ satisfying $d\rho(z) \neq 0$ if $\rho(z) =0$.\\

For $0 \leq q \leq n \ (\geq 2)$, let $C^{\infty}_q({\bar \Omega})$ be the vector space of $(0,q)$-forms with coefficients in  $C^{\infty}({\bar \Omega})$, the vector space of complex-valued functions $f$ such that $f$ is infinitely differentiable on an open neighborhood $U_f$ of ${\bar \Omega}$.   The Cauchy-Riemann operator $\bar \partial$ gives rise to (a special version of) the {\it Dolbeault complex} (or the {\it Cauchy-Riemann complex})
$$ 0 \rightarrow  C^{\infty}_0({\bar \Omega})\overset {{\bar \partial_0}} \rightarrow  C^{\infty}_1({\bar \Omega})\rightarrow\ldots \overset {{\bar \partial}_{n-1}} \rightarrow   C^{\infty}_n({\bar \Omega})\rightarrow 0.$$
Using the normalized volumetric measure $\nu$ on ${\bar \Omega}$ one can define an inner product on $C^{\infty}_q({\bar \Omega})$ in a natural way (refer to \cite[Chapter 2, Section 2.1]{U}).
Let $L^2_{q}(\Omega)$ be the Hilbert space completion of $C^{\infty}_q({\bar \Omega})$ in this inner product, with the corresponding norm on $L^2_{q}(\Omega)$ being denoted by $\|\cdot\|$ (for any $q$).
The closure of $\bar \partial_q$ will still be denoted by $\bar \partial_q$; thus $\bar \partial_q$ is a densely defined closed (linear) operator from $L^2_{q}(\Omega)$ into $L^2_{q+1}(\Omega)$. The Hilbert space adjoint of ${\bar \partial}_q$ will be denoted by ${\bar \partial}^*_{q+1}$ (unlike ${\bar \partial}^*_{q}$ in (2.1.13) of \cite{U} which, in view of the subsequent formulas employed there, is a notational inaccuracy). 
The ($q$th) {\it ${\bar \partial}$-Neumann Laplacian} is defined by ${\square}_q = {\bar \partial}_{q-1}{\bar \partial}^*_{q}+ {\bar \partial}^*_{q+1}{\bar \partial}_{q}$ (with $\bar \partial_n$, $\bar \partial_{-1}$, $\bar \partial^*_{n+1}$ and $\bar \partial^*_0$ being interpreted as zero operators). For $1 \leq q \leq n$, ${\square}_q$ turns out to be invertible with a bounded inverse $N_q$ (refer to \cite{F-Ko}, \cite{Ho}); the operator $N_q$ is referred to as the ($q$th) {\it $\bar \partial$-Neumann operator}.\\

For $0 \leq q \leq n \ (\geq 2)$, let $R^{\infty}_q(\partial \Omega)$ be the vector space obtained by restricting the members of $C^{\infty}_q({\bar \Omega})$ to ${\partial \Omega}$. If  $f = \sum_{i_1<\cdots<i_q} \phi_{i_1,\ldots,i_q}z_{i_1}\wedge\cdots\wedge z_{i_q}$ and  $g = \sum_{i_1<\cdots<i_q} \psi_{i_1,\ldots,i_q}z_{i_1}\wedge\cdots\wedge z_{i_q}$ (in the standard notation) are in $R^{\infty}_q(\partial \Omega)$, then $f$ is said to be {\it pointwise orthogonal} to $g$ if $\sum_{i_1<\cdots<i_q} \phi_{i_1,\ldots,i_q}(b){\overline {\psi_{i_1,\ldots,i_q}(b)}} =0$ for every $b \in \partial \Omega$ (notation: $f \perp g$). If $N^{\infty}_q(\partial \Omega)$ is the vector space $\{ f \in R^{\infty}_q(\partial \Omega): f \wedge ({\bar \partial} \rho|\partial \Omega) = 0\}$, then we declare $C^{\infty}_q({\partial \Omega})$ to be the vector space $\{f \in R^{\infty}_q(\partial \Omega): f \perp g {\rm \ for\ all\ } g \in N^{\infty}_q(\partial \Omega)\}$; it is to be noted that $C^{\infty}_n({\partial \Omega}) = \{0\}$.  The Cauchy-Riemann operator $\bar \partial$ induces the tangential Cauchy-Riemann operator ${\bar \partial}_b$ (refer to \cite{F-Ko} and \cite{Ko-R}) that gives rise to (a special version of) the {\it Kohn-Rossi complex} (or the {\it tangential Cauchy-Riemann complex})  
$$ 0 \rightarrow  C^{\infty}_0({\partial \Omega})\overset {{\bar \partial}_{b,0}} \rightarrow  C^{\infty}_1({\partial \Omega})\rightarrow\ldots \overset {{\bar \partial}_{b,n-2}} \rightarrow   C^{\infty}_{n-1}({\partial \Omega})\rightarrow 0.$$
The vector space $C^{\infty}_q({\partial \Omega})$ can be equipped naturally with an inner product by using the normalized surface area measure $\sigma$ on $\partial \Omega$ (refer to \cite[Chapter 2, Section 2.2]{U}).
Let $L^2_{q}(\partial \Omega)$ be the Hilbert space completion of $C^{\infty}_q({\partial \Omega})$ in this inner product.
The closure of $\bar \partial_{b,q}$ will still be denoted by $\bar \partial_{b,q}$; thus $\bar \partial_{b,q}$ is a densely defined closed (linear) operator from $L^2_{q}(\partial \Omega)$ into $L^2_{q+1}(\partial \Omega)$. The Hilbert space adjoint of ${\bar \partial}_{b,q}$ will be denoted by ${\bar \partial}^*_{b,q+1}$ (with the notational inaccuracy in (2.2.9) of \cite{U} noted). 
The ($q$th) {\it Kohn Laplacian} is defined by ${\square}_{b,q} = {\bar \partial}_{b,q-1}{\bar \partial}^*_{b,q}+ {\bar \partial}^*_{b,q+1}{\bar \partial}_{b,q}$ (with $\bar \partial_{b,n-1}$, $\bar \partial_{b,-1}$, $\bar \partial^*_{b,n}$ and $\bar \partial^*_{b,0}$ being interpreted as zero operators).  For $1 \leq q \leq n-1$, ${\square}_{b,q}$ turns out to be invertible with a bounded inverse $N_{b,q}$ (refer to \cite{F-Ko}, \cite{Ko}); the operator $N_{b,q}$ is referred to as the ($q$th) {\it  complex Green operator} or the ($q$th) {\it tangential Neumann operator}.\\

Let $W^{-1}_1(\Omega)$ be the vector space of $(0,1)$-forms $f$ with coefficients in the Sobolev space $W^{-1}(\Omega)$ of order $-1$, and let $\|f\|^2_{-1}$ be the sum of squares of the $W^{-1}(\Omega)$ norms of the coefficients of $f$. One says that a {\it compactness estimate}  holds (for $\Omega$) if for every positive $\epsilon$ there exists $C(\epsilon)$ such that 
$$ \|f\|^2 \leq \epsilon \{\|\bar \partial_1f\|^2 + \|\bar \partial^*_1f\|^2\} + C(\epsilon) \|f\|^2_{-1}$$
for all $(0,1)$-forms $f$ that lie in ${\rm Domain}(\bar \partial_1) \cap {\rm Domain}(\bar \partial^*_1)$ $(\subset L^2_1(\Omega) \subset W^{-1}_1(\Omega))$.\\

One says that $\partial \Omega$ satisfies the {\it Catlin property $(P)$} if for every positive $M$ there exists a plurisubharmonic function $\lambda$ in  $C^{\infty}({\bar \Omega})$ with $0 \leq \lambda \leq 1$ such that 
$$ \sum_{j,k =1}^n\frac{\partial^2\lambda}{\partial z_j \partial {\bar z}_k}(b) t_j {\bar t}_k \geq M\{|t_1|^2+\cdots +|t_n|^2\}$$
for all points $t = (t_1,\ldots,t_n)$ in $\C^n$ and for all points $b$ of $\partial \Omega$.\\

{\bf Remark 3.1}. If a bounded pseudoconvex domain $\Omega$ has real analytic boundary $\partial \Omega$, then it follows from \cite[Lemma 2]{Die-F} and \cite[Theorem 2]{Cat} that $\partial \Omega$ satisfies the Catlin property $(P)$; in particular, $\partial \Omega_p$ satisfies the Catlin property $(P)$.\\

The closure of $A(\Omega_p)$ in $L^2(\nu_p)$, where $\nu_p$ is the normalized volumetric measure on $\bar \Omega_p$, will be referred to as the {\it Bergman space} of $\Omega_p$ and will be denoted by $A^2(\nu_p)$. The tuple of multiplications by the coordinate functions $z_i$ on $A^2(\nu_p)$ will be denoted by $M_{\nu_p,z}$. Let $\tilde P_{\nu_p}$ be the orthogonal projection of $\L^2(\nu_p)$ onto $A^2(\nu_p)$ and let, for $\phi \in L^{\infty}(\nu_p)$, $\tilde N_{\nu_p,\phi}$ denote the operator of multiplication by $\phi$ on $L^2(\nu_p)$. We let $\tilde T_{\nu_p,\phi}$ stand for ${\tilde P_{\nu_p}}{\tilde N}_{\nu_p,\phi}|A^2(\nu_p)$ and refer to $\tilde T_{\nu_p,\phi}$ as a {\it Bergman-Toeplitz operator}. The adjoint of the Bergman-Toeplitz operator $\tilde T_{\nu_p,\phi}$ (resp.  $M_{\mu_p,z}$-Toeplitz operator $T_{\mu_p,\phi}$ of Remark 2.6) equals $\tilde T_{\nu_p,\bar \phi}$ (resp. $T_{\mu_p,\bar \phi}$).\\

{\bf Remark 3.2}. The domain $\Omega_p$ is starlike with respect to the origin, and any $f \in A(\Omega_p)$ can be approximated uniformly on $\bar \Omega_p$ by the sequence $\{f_m\}$ of functions $f_m$ in $\O({\overline {\Omega}_p})$ where $f_m(z) = f((1-\frac{1}{m})z)$. Further,  $\bar \Omega_p$ is polynomially convex so that any function such as $f_m$ that is holomorphic on an open neighborhood of $\bar \Omega_p$ is the uniform limit (on $\bar \Omega_p$) of polynomials by the Oka-Weil approximation theorem (see \cite[Chapter VI, Theorem 1.5]{Ran}. It is then clear that the Hardy space $H^2(\sigma_p)$ (resp. Bergman space $A^2(\nu_p)$) as defined previously is really the closure of polynomials in $L^2(\sigma_p)$ (resp. $L^2(\nu_p)$) with the constant function $1$ in $H^2(\sigma_p)$ (resp. $A^2(\nu_p)$) being a cyclic vector for $M_{\sigma_p,z}$ (resp. $M_{\nu_p,z}$). The multiplication tuple $M_{\sigma_p,z}$ (resp. $M_{\nu_p,z}$) can be looked upon as a multivariable weighted shift, with the positive weights of $M_{\sigma_p,z}$ (resp. $M_{\nu_p,z}$) being computed by checking the action of each $M_{\sigma_p,z_i}$ (resp. $M_{\nu_p,z_i}$) on the members of the orthonormal basis obtained by applying the Gram-Schmidt process to the constant function $1$ and the powers of $z_i$ in the Hardy space $H^2(\sigma_p)$ (resp. Bergman space $A^2(\nu_p)$)(refer to \cite{Je-L}). For an arbitrary $\Omega_p$, such computations can turn out to be formidable as can be gathered, for example, by referring to similar computations carried out in \cite{Cr-R} in the context of `complex ellipsoids' in $\C^n$.  \\

{\bf Proposition 3.3}. The semicommutator $\tilde T_{\nu_p,\phi}\tilde T_{\nu_p,\psi}- \tilde T_{\nu_p,\phi\psi}$ of the Bergman-Toeplitz operators $\tilde T_{\nu_p,\phi}$ and $\tilde T_{\nu_p,\psi}$ is compact for any continuous functions $\phi$ and $\psi$ on $\bar \Omega_p$.
\begin{proof} In view of Remark 3.1, $\partial \Omega_p$ satisfies the Catlin property $(P)$. The Catlin property $(P)$ implies that a compactness estimate holds for $\Omega_p$ (refer to \cite[Theorem 1]{Cat}). That in turn implies that the $\bar \partial$-Neumann operator $N_1$ correponding to $\Omega_p$ is compact (refer to \cite[Lemma 11]{Bo-S}). Now arguing exactly as in \cite[Lemma 2.1.24]{U}, one proves that $\bar \partial^*_1N_1$ is a compact operator. (The symbol $\bar \partial^*_0$ in the proof of \cite[Lemma 2.1.24]{U} should be corrected to $\bar \partial^*_2$). Next, using the compactness of $\bar \partial^*_1N_1$ and arguing as in \cite[Lemma 2.1.22]{U} and \cite[Theorem 4.1.18]{U}, one proves that $(I-{\tilde P}_{\nu_p}){\tilde N}_{\nu_p,\phi}|A^2(\nu_p) : A^2(\nu_p) \rightarrow L^2(\nu_p)$ is compact for any  $\phi$ that is continuous on $\bar \Omega_p$. And, as in \cite[Corollary 4.1.21]{U}, that leads to the compactness of  the semicommutator $\tilde T_{\nu_p,\phi}\tilde T_{\nu_p,\psi}- \tilde T_{\nu_p,\phi\psi}$ for any continuous functions $\phi$ and $\psi$ on $\bar \Omega_p$ .
\end{proof}
{\bf Corollary 3.4}. The commutator $\tilde T_{\nu_p,\phi}\tilde T_{\nu_p,\psi}- \tilde T_{\nu_p,\psi}\tilde T_{\nu_p,\phi}$ of the Bergman-Toeplitz operators $\tilde T_{\nu_p,\phi}$ and $\tilde T_{\nu_p,\psi}$ is compact for any continuous functions $\phi$ and $\psi$ on $\bar \Omega_p$; in particular, the multiplication tuple $M_{\nu_p,z}$ is essentially normal.\\

{\bf Proposition 3.5}. Let $n \geq 3$. For $\Omega_p \subset \C^n$, the semicommutator $T_{\sigma_p,\phi}T_{\sigma_p,\psi}- T_{\sigma_p,\phi\psi}$ of the $M_{\sigma_p,z}$-Toeplitz operators $T_{\sigma_p,\phi}$ and $T_{\sigma_p,\psi}$ is compact for any continuous functions $\phi$ and $\psi$ on $\partial \Omega_p$.
\begin{proof} In view of Remark 3.1, $\partial \Omega_p$ satisfies the Catlin property $(P)$. It follows from \cite[Theorem 1.4]{Ra-St} that the tangential Neumann operator $N_{b,1}$ corresponding to $\Omega_p$ is compact. Now arguing exactly as in the Bergman case, one proves an analog of \cite[Lemma 2.2.19]{U} to obtain that $\bar \partial^*_{b,1}N_{b,1}$ is a compact operator. Next, using the compactness of $\bar \partial^*_{b,1}N_{b,1}$ (and arguing as in the Bergman case) one establishes analogs of \cite[Lemma 2.2.18]{U} and \cite[Theorem 4.2.17]{U} to obtain that $(I-P_{\sigma_p}){N}_{\sigma_p,\phi}|H^2(\sigma_p) : H^2(\sigma_p) \rightarrow L^2(\sigma_p)$ is compact for any  $\phi$ that is continuous on $\partial \Omega_p$. And that leads to an analog of \cite[Corollary 4.2.20]{U}, yielding the compactness of  the semicommutator
$T_{\sigma_p,\phi}T_{\sigma_p,\psi}- T_{\sigma_p,\phi\psi}$ for any continuous functions $\phi$ and $\psi$ on $\partial \Omega_p$.
\end{proof}
{\bf Corollary 3.6}. Let $n \geq 3$. For $\Omega_p \subset \C^n$, the commutator $T_{\sigma_p,\phi}T_{\sigma_p,\psi}- T_{\sigma_p,\psi}T_{\sigma_p,\phi}$ of the $M_{\sigma_p,z}$-Toeplitz operators $T_{\sigma_p,\phi}$ and $T_{\sigma_p,\psi}$ is compact for any continuous functions $\phi$ and $\psi$ on $\partial \Omega_p$; in particular, the multiplication tuple $M_{\sigma_p,z}$ is essentially normal.\\

The author does not know whether the tangential Neumann operator $N_{b,1}$ corresponding to an arbitrary $\Omega_p \subset \C^2$ is compact; as such a different strategy is adopted below to prove the essential normality of the multiplication pair $M_{\sigma_p,z} \in ({\mathcal B}(H^2(\sigma_p))^2$ for any $\Omega_p \subset \C^2$.\\

{\bf Proposition 3.7}. For $\Omega_p \subset \C^2$, the multiplication pair $M_{\sigma_p,z}$ is essentially normal.
\begin{proof} Since $\Omega_p\ (\subset \C^2)$ is a pseudoconvex complete Reinhardt domain with real analytic boundary, it follows from the work of Sheu in \cite{Sh} that there is a $*$-isomorphism $\Psi$ of the  C$^*$-algebra $\A$ generated by the set $\{{\tilde T}_{\nu_p,\phi}: \phi \ {\rm is \ continuous \ on\ } {\bar \Omega_p}\}$ with the C$^*$-algebra $\B$ generated by the set $\{T_{\sigma_p,\phi}: \phi \ {\rm is \ continuous \ on\ } {\partial \Omega_p}\}$. In view of Remark 3.2 and \cite[Corollary 13]{Je-L}, the  C$^*$-algebras $\A$ and $\B$ are irreducible. Let ${\mathcal K}(A^2(\nu_p))$ (resp. ${\mathcal K}(H^2(\sigma_p))$) be the C$^*$-algebra of compact operators on $A^2(\nu_p)$ (resp. $H^2(\sigma_p)$). As $\A$ has (by Corollary 3.4) a non-trivial intersection with ${\mathcal K}(A^2(\nu_p))$ and as $\A$ is irreducible, $\A$ contains ${\mathcal K}(A^2(\nu_p))$ (refer to \cite{Co2}). Consider $\Psi|{\mathcal K}(A^2(\nu_p)): {\mathcal K}(A^2(\nu_p))\rightarrow {\mathcal B}(H^2(\sigma_p))$. Since $\Psi({\mathcal K}(A^2(\nu_p)))$ is an ideal of $\B$ and since $\B$ is irreducible, $(\Psi|{\mathcal K}(A^2(\nu_p),H^2(\sigma_p))$ is an irreducible representation of ${\mathcal K}(A^2(\nu_p))$. It follows then from \cite[ Corollary 16.12]{Co2} that $\Psi({\mathcal K}(A^2(\nu_p)) = {\mathcal K}(H^2(\sigma_p))$.  Letting $T_i = \Psi^{-1}(M_{\sigma_p,z_i})$, it is clear that the compactness of $M_{\sigma_p,z_i}^*M_{\sigma_p,z_j}- M_{\sigma_p,z_j}M_{\sigma_p,z_i}^* \in {\mathcal B}(H^2(\sigma_p))$ would follow from that of $T_i^*T_j-T_jT_i^* \in {\mathcal B}(A^2(\nu_p))$. However, the compactness of $T_i^*T_j-T_jT_i^* \in {\mathcal B}(A^2(\nu_p))$ can be deduced easily from the result of Proposition 3.3 and the fact that the uniform limit of compact operators is compact. 
\end{proof}

The results of Corollary 3.6 and Proposition 3.7 allow us to bring all the results in \cite{Di-E-Ev} and \cite{E-Ev} related to an essentially normal regular $A$-isometry to bear upon the multiplication tuple $M_{\sigma_p,z}$; we highlight in Remark 3.6 below a couple of implications of the results in those works.\\

{\bf Remark 3.8}. (a) Let $\T_a\left(M_{\sigma_p, z}\right)$ be the set $\{T_{\sigma_p,\phi}:\phi \in H^{\infty}_{A(\Omega_p)}(\sigma_p)\}$. For $\phi \in L^{\infty}$, let $H_{\sigma_p,\phi}$ be the {\it Hankel operator} from $H^2(\sigma_p)$ to $H^2(\sigma_p)^{\perp} =L^2(\sigma_p) \ominus H^2(\sigma_p)$ defined by  $H_{\sigma_p,\phi} = (I-P_{\sigma_p})N_{\sigma_p,\phi}|H^2(\sigma_p)$. In view of the observations in the proof of \cite[Corollary 3.3]{Di-E-Ev} and in view of \cite[Corollary 5.1]{E-Ev}, one has that an operator $S \in  
{\mathcal B}(H^2(\sigma_p))$ is in the essential commutant of $\T_a\left(M_{\sigma_p, z}\right)$ if and only if $S$ equals $T_{\sigma_p,\phi} +K$ for some compact operator $K$ on 
$H^2(\sigma_p)$ and some $\phi$ in $L^{\infty}(\sigma_p)$ for which the Hankel operator $H_{\sigma_p,\phi}$ is compact.\\
(b) From \cite[Proposition 3.10]{Di-E-Ev} one can deduce the existence of a short exact sequence of C$^*$-algebras
$$ 0 \rightarrow {\mathcal K}(H^2(\sigma_p)) \overset {\iota}\rightarrow {\B} \overset{\pi}\rightarrow C(\partial \Omega_p) \rightarrow 0$$
where ${\mathcal K}(H^2(\sigma_p))$ and $\B$ are as in the proof of Proposition 3.7, $\iota$ is the inclusion map, and $\pi$ is a unital $*$-homomorphism satisfying $\pi(T_{\sigma_p,\phi}) = \phi$ for any $\phi \in C(\partial \Omega_p)$. \\

\section{$\partial\Sigma_p$-isometries}
Let $p= (p_1,p_2,\ldots,p_n)$ be an $n$-tuple of $m_i$-tuples $p_i= (p_{i,1},p_{i,2},\ldots,p_{i,m_i})$ where $p_{i,1}$,..., $p_{i,m_i}$ (with $m_i \geq 1$) are positive integers. The subset $\Sigma_p$ of $\C^n$ is defined by $\Sigma_p = \{z =(z_1,z_2,\ldots,z_n) \in \C^n: 
\sum_{i=1}^n\sum_{j=1}^{m_i} |z_i|^{2p_{i,j}} < 1\}$. The set $\Sigma_p$ is easily seen to be a convex complete Reinhardt domain in $\C^n$ with the real analytic boundary $\partial\Sigma_p=\{z =(z_1,z_2,\ldots,z_n) \in \C^n: \sum_{i=1}^n\sum_{j=1}^{m_i} |z_i|^{2p_{i,j}} =1\}$. We use the symbol $\Sigma^{(n)}$ to denote the class of domains $\Sigma_p$ in $\C^n$. Trivially, $\Sigma^{(n)}$ is a superclass of the class $\Omega^{(n)}$. The domain $\Sigma_p$ is a so-called {\it complex ellipsoid} in case $m_i =1$ for each $i$; we also note that, for any $n$, $\B_n \in \Sigma^{(n)} \setminus \Omega^{(n)}$.\\

{\bf Definition 4.1}.  If $S = (S_1,\ldots,S_n)$ is a subnormal $n$-tuple of (commuting) operators $S_i$ in ${\mathcal B}({\mathcal H})$ such that the spectral measure $\rho_N$ of the minimal normal extension $N$ of $S$ is supported on $\partial\Sigma_p$, then $S$ is called a {\it $\partial\Sigma_p$-isometry}. \\

{\bf Remark 4.2}. The statements (and proofs) of Propositions 2.4, 3.3, 3.5, 3.7 along with those of Corollaries 3.4, 3.6 hold and the contents of Remarks 3.1, 3.2 remain applicable with $\Sigma_p$ in place of $\Omega_p$ and with the obvious corresponding interpretations of $A(\Sigma_p)$, $\sigma_p$, $\nu_p$, $H^2(\sigma_p)$, $A^2(\nu_p)$, $M_{\sigma_p,z}$, $M_{\nu_p,z}$, $\A$ and $\B$. We also note that any $\Sigma_p$-isometry $S \in {\mathcal B}({\mathcal H})^n$ is an $A(\Sigma_p)$-isometry. (Indeed, if $\mu_p$ is a scalar spectral measure of the minimal normal extension $N$ of $S$, then $\mu_p$ is supported on $\partial\Sigma_p$ where $\partial\Sigma_p$ is the Shilov boundary of $A(\Sigma_p)$ by the analog of Proposition 2.4 for $\Sigma_p$; thus we need only check that $A(\Sigma_p)$ is contained in the restriction algebra ${\mathcal R}_S$ of $S$. Let $f \in A(\Sigma_p)$. Choosing $f_m$ as in Remark 3.2 and using the Taylor functional calculus for $S$ (refer to \cite{T}), one has that $f_m(N)|{\mathcal H}=f_m(S)\in {\mathcal B}({\mathcal H})$. Since the sequence $\{f_m\}$ converges to $f$ uniformly on $\partial \Sigma_p$, it is clear that $f(N){\mathcal H} \equiv \Psi_N(f){\mathcal H}$ is contained in ${\mathcal H}$). Thus $\Sigma_p$-isometries, like the less general $\Omega_p$-isometries, are examples of essentially normal $A$-isometries, but $\Omega_p$-isometries come with an added bonus of regularity.\\

{\bf Remark 4.3}. The weak$^*$ closure $H^{\infty}_{A(\Sigma_p)}(\sigma_p)$ of $A(\Sigma_p)$ in $L^{\infty}(\sigma_p)$ can be identified with the algebra $H^{\infty}(\sigma_p)$ of the non-tangential boundary limits of the members of $H^{\infty}(\Sigma_p)$ where $H^{\infty}(\Sigma_p)$ is the algebra of bounded holomorphic functions on $\Sigma_p$. Indeed, $H^{\infty}(\Sigma_p)$ can be shown to be a weak*-closed subalgebra of $L^{\infty}(\sigma_p)$ and the map ${\tilde r}_{\sigma_p} : H^{\infty}(\Sigma_p) \rightarrow L^{\infty}(\sigma_p)$ that associates with any $f \in H^{\infty}(\Sigma_p)$ its non-tangential boundary limit can be shown to be an isometric and a  weak$^*$-continuous algebra homomorphism  as in the argument provided in the discussion preceding \cite[Corollary 4.8]{Di-E-Ev}; further, also as per  the argument there, the inclusion $H^{\infty}_{A(\Sigma_p)}(\sigma_p) \subset H^{\infty}(\sigma_p) (= {\tilde r}_{\sigma_p}(H^{\infty}(\Sigma_p))$ holds. For the other way inclusion one uses that ${\tilde r}_{\sigma_p}$ is weak$^*$-continuous and that any function $f \in  H^{\infty}(\Sigma_p)$ can be approximated in the weak$^*$ topology of $L^{\infty}(\sigma_p)$ by the sequence $\{f_m\}$ where $f_m$ are as in Remark 3.2. \\

An intrinsic characterization of $\partial \Sigma_p$-isometries can be provided using the results of \cite{At-P}. If $q(z,w)= \sum_{\alpha,\beta} a_{\alpha,\beta}z^{\alpha}w^{\beta}$ is a polynomial in the variables $z=(z_1,\ldots,z_n)$ and $w=(w_1,\ldots,w_n)$ with real coefiicients $a_{\alpha,\beta}$, then for any $n$-tuple $T = (T_1,\ldots,T_n)$ of commuting operators in ${\mathcal B}({\mathcal H})$ we interpret $(q(z,w))(T,T^*)$ to be the operator $\sum_{\alpha,\beta} a_{\alpha,\beta}T^{*\beta}T^{\alpha}$. Since the Taylor spectrum of the minimal normal extension of  a $\partial \Sigma_p$-isometry $S$ is contained in $\partial \Sigma_p$, it follows by a result of Curto \cite{Cu} that the Taylor spectrum of $S$ is contained in the polynomial convex hull of $\partial \Sigma_p$, which is the closure $\bar {\Sigma}_p$ of $\Sigma_p$. As  $\bar {\Sigma}_p$ is contained in the closed unit polydisk in $\C^n$ centered at the origin, the spectral projection property of the Taylor spectrum implies that any coordinate $S_i$ of $S$ has its spectrum contained in the unit disk in $\C$ centered at the origin so that the spectral radius $r_{S_i}$ of $S_i$ cannot exceed $1$. Since $S_i$ is subnormal, the norm of $S_i$ must equal $r_{S_i}$ (refer to \cite{Co1}) and hence $S_i$ is  a contraction. The following result is now a consequence of \cite[Proposition 7]{At-P} and the observations in the proof of \cite[Proposition 8]{At-P}.\\ 

{\bf Proposition 4.4}. Let $S = (S_1,\ldots,S_n)$ be an $n$-tuple of commuting operators $S_i$ in ${\mathcal B}({\mathcal H})$. The statements (i) and (ii) below are equivalent:\\
(i) $S$ is a $\partial \Sigma_p$-isometry.\\
(ii) (a) $(\Pi_{i=1}^n [1-z_iw_i]^{k_i})(S,S^*) \geq 0$ for all integers $k_i \geq 0$.\\
(b) $(1-\sum_{i=1}^n\sum_{j=1}^{m_i} z_i^{p_{i,j}}w_i^{p_{i,j}})(S,S^*) = 0$.\\

The condition (b) in part (ii) of Proposition 4.4 can simply be written as 
$$ I-\sum_{i=1}^n\sum_{j=1}^{m_i} S_i^{*p_{i,j}}S_i^{p_{i,j}} = 0$$
and, as shown below, by itself characterizes  a $\partial \Sigma_p$-isometry for a special type of $\Sigma_p$. We consider those $\Sigma_p$ (with $p = (p_1,\ldots,p_n)$) for which each $p_i$ has at least one integer coordinate equal to $1$; we use the symbol ${\tilde \Sigma}_p$ to denote any such $\Sigma_p$ and note that ${\tilde \Sigma}_p$ is strictly pseudoconvex. The unit ball 
${\mathbb B}_n = \Sigma_{(p_1,\ldots,p_n)}$ 
with $p_i = (1)$ for every $i$ is a special example of such a domain; we note that $\partial {\mathbb B}_n$-isometries are precisely spherical isometries. The next proposition provides a characterization of a $\partial {\tilde \Sigma}_p$-isometry that is a generalization of that of a spherical isometry.\\

{\bf Proposition 4.5}. Let $S = (S_1,\ldots,S_n)$ be an $n$-tuple of commuting operators $S_i$ in ${\mathcal B}({\mathcal H})$. The statements (i) and (ii) below are equivalent:\\
(i) $S$ is a $\partial {\tilde \Sigma}_p$-isometry.\\
(ii) $(1-\sum_{i=1}^n\sum_{j=1}^{m_i} z_i^{p_{i,j}}w_i^{p_{i,j}})(S,S^*) = 0$.
\begin{proof}
The implication $(i) \implies (ii)$ is trivial. To prove $(ii) \implies (i)$, we need only show that the condition  $(1-\sum_{i=1}^n\sum_{j=1}^{m_i} z_i^{p_{i,j}}w_i^{p_{i,j}})(S,S^*) = 0$ guarantees the condition (ii) (a) of Proposition 4.4, viz, $(\Pi_{i=1}^n[1-z_iw_i]^{k_i})(S,S^*) \geq 0$ for all integers $k_i \geq 0$. We assume without any loss of generality that $p_{i,1} = 1$ for each $i$. Let $q_i(z,w) = \sum_{j=2}^{m_i}z_i^{p_{i,j}}w_i^{p_{i,j}}+ \sum_{k\neq i}\sum_{j=1}^{m_k}z_k^{p_{k,j}}w_k^{p_{k,j}}$. That the condition (ii) (a) of Proposition 4.4 holds follows by observing that $(\Pi_{i=1}^n[1-z_iw_i]^{k_i})(S,S^*)$ can be written as $(\Pi_{i=1}^n[\{1-\sum_{i=1}^n\sum_{j=1}^{m_i} z_i^{p_{i,j}}w_i^{p_{i,j}}\}+q_i(z,w)]^{k_i})(S,S^*)$.
\end{proof}

We now turn to examining the intertwining of two $\partial \Omega_p$-isometries. By choosing $S=T$ in the following proposition, one obtains a commutant lifting theorem for a $\partial \Omega_p$-isometry.\\

{\bf Proposition 4.6}. Let $S =(S_1,\ldots,S_n) \in {\mathcal B}({\mathcal H})^n$ and $T =(T_1,\ldots,T_n) \in {\mathcal B}({\mathcal K})^n$ be $\partial \Omega_p$-isometries, and let 
$M =(M_1,\ldots,M_n) \in {\mathcal B}({\tilde {\mathcal H}})^n$ and $N =(N_1,\ldots,N_n) \in {\mathcal B}({\tilde {\mathcal K}})^n$ be the minimal normal extensions of $S$ and $T$,
respectively. If $X: H \rightarrow K$ is a bounded linear map intertwining $S$ and $T$ so that $XS_i = T_iX$ for all $i$, then $X$ lifts to a bounded linear map $\tilde X: {\tilde H} \rightarrow {\tilde K}$ intertwining $M$ and $N$ and satisfying $\|{\tilde X}\| =\|X\|$.
\begin{proof}
Since the Taylor spectra of $M$ and $N$ are contained in $\partial \Omega_p$, by the result of Curto \cite{Cu} (mentioned earlier) the Taylor spectra of $S$ and $T$ are contained in the polynomial convex hull of $\partial \Omega_p$, which is $\bar \Omega_p$. Let $f \in A(\Omega_p)$. For any positive integer $m \geq 2$, $f_m$ defined by $f_m(z) = f((1-\frac{1}{m})z)$ is holomorphic on an open neighborhood of $\bar \Omega_p$. If $X$ intertwines $S$ and $T$, then it follows by the Taylor functional calculus (see \cite[Proposition 4.5]{T}) that $Xf_m(S) = f_m(T)X$. If $\rho_M$ (resp. $\rho_N$) is the spectral measure of $M$ (resp. $N$), then $\rho_S= P_{\mathcal H}\rho_M|{\mathcal H}$ (resp. $\rho_T= P_{\mathcal K}\rho_N|{\mathcal K})$ is the semi-spectral measure of $S$ (resp. $T$) with $P_{\mathcal H}$ and $P_{\mathcal K}$ being appropriate projections, and for any $u \in \mathcal H$ and any $v \in \mathcal K$ one has
$$ \|f_m(S)u\|^2 = \int|f_m(z)|^2 d\langle\rho_S(z)u,u\rangle$$
and
$$ \|f_m(T)v\|^2 = \int|f_m(z)|^2 d\langle\rho_T(z)v,v\rangle.$$
Letting $v = Xu$ and using $Xf_m(S) = f_m(T)X$, one obtains
$$ \int|f_m(z)|^2 d\langle\rho_T(z)Xu,Xu\rangle \leq \|X\|^2\int|f_m(z)|^2 d\langle\rho_S(z)u,u\rangle,$$
which, upon letting $m$ tend to infinity, yields
$$ \int|f(z)|^2 d\langle\rho_T(z)Xu,Xu\rangle \leq \|X\|^2\int|f(z)|^2 d\langle\rho_S(z)u,u\rangle.$$
Consider $\eta(\cdot) = \langle\rho_T(\cdot)Xu,Xu\rangle + \langle\rho_S(\cdot)u,u\rangle$. One has by Proposition 2.5 that $(A(\Omega_p)|{\partial \Omega_p}, {\partial \Omega_p}, \eta)$ is a regular triple. Thus if $\phi$ is any positive continuous function on $\partial \Omega_p$, then there exists a sequence of functions $\{\phi_m\}_{m \geq 1}$ in $A(\Omega_p)$ such that $|\phi_m| < \sqrt \phi$ on $\partial \Omega_p$ and $\lim_{m \rightarrow \infty} |\phi_m| = \sqrt \phi$ $\eta$-almost everywhere. Replacing $f$ by $\phi_m$ in the last integral inequality and letting $m$ tend to infinity, one obtains
$$ \int\phi(z) d\langle\rho_T(z)Xu,Xu\rangle \leq \|X\|^2\int\phi(z) d\langle\rho_S(z)u,u\rangle.$$
That yields $\langle\rho_T(\cdot)Xu,Xu\rangle \leq \|X\|^2\langle\rho_S(\cdot)u,u\rangle$ for every $u$ in $\mathcal H$. The desired conclusion now follows by appealing to \cite[Lemma 4.1]{M}.
\end{proof}

{\bf Remark 4.7}. Requiring $X$ to be of a special type in Proposition 4.6 may guarantee the lift $\tilde X$ of $X$ also to be of that special type.  Indeed, arguing as in \cite[Theorem 5.2]{M}, one can establish the following facts: If $X$ is isometric, then so is ${\tilde X}$; if $X$ has dense range, then so has ${\tilde X}$; if $X$ is bijective, then so is ${\tilde X}$. If a bounded linear map $X$ that intertwines $n$-tuples $S$ and $T$ is invertible (resp. unitary), then we refer to $S$ and $T$ as being {\it similar} (resp. {\it unitarily equivalent}). It follows from \cite[Lemma 1]{At2} and Proposition 4.6 above that if $\partial\Omega_p$-isometries $S$ and $T$  are intertwined by a  bounded linear map $X$ that is injective and has dense range (that is, if $S$ and $T$ are {\it quasisimilar}), then the minimal normal extensions of $S$ and $T$ are unitarily equivalent (cf. \cite[Proposition 9]{At2}).\\

In the light of Remark 3.2, it is natural to investigate analogs of Proposition 4.6 for a pair of subnormal tuples, one of which is a cyclic $\partial \Omega_p$-isometry. It is a standard fact of the subnormal operator theory (refer, for example, to \cite{Ha}) that any cyclic subnormal tuple $S$ is (up to unitary equivalence) a multiplication tuple $M_{\theta,z}$ on the closure $P^2(\theta)$ of polynomials in $L^2(\theta)$ for some compactly supported positive regular Borel measure $\theta$; in case $S$ happens to be a cyclic $\partial\Omega_p$-isometry, $\theta$ must be supported on $\partial \Omega_p$.\\  

{\bf Hereafter, $T=M_{\theta,z}$  stands for a fixed cyclic $\partial\Omega_p$-isometry with $\theta$ supported on $\partial\Omega_p$ and having no atoms on $\partial\Omega_p$}.\\

To discuss subnormal tuples $S$ quasisimilar to $T=M_{\theta,z}$, we need only consider $S = M_{\eta,z}$ for some compactly supported positive regular Borel measure $\eta$ on $\C^n$ as is justified by \cite[Proposition 1]{At1}.\\

Arguing almost verbatim along the lines of \cite[Section 4]{At3} (refer also to \cite{At1}), where the context was that of  strictly pseudoconvex domains, one can establish Lemmas 4.8 and 4.9 and Propositions 4.10 and 4.11 below. That one can use polynomials in the statements of those lemmas and propositions is a pleasant consequence of our observations in Remark 3.2.  We point out that, like in the proof of \cite[Lemma 4.5]{At3}, one needs to appeal in the proof of Lemma 4.9 below to \cite[Corollary 2.8]{Di-E}, which is a consequence of some refinement in \cite{Di-E} of Aleksandrov's work in \cite{A}; the requirement that $\theta$ have no atoms on $\partial \Omega_p$ stems from the need to apply \cite[Corollary 2.8]{Di-E}.\\

{\bf Lemma 4.8}.
Let $S$ be a cyclic subnormal tuple so that $S$ can be identified with $M_{\eta,z}$ for some compactly supported positive regular Borel measure $\eta$ on $\C^n$. If there exists a bounded linear map $Y: {P}^2(\theta) \rightarrow {P}^2(\eta)$ with dense range such that $YM_{\theta,z} = M_{\eta,z}Y$, then there exists a cyclic vector $g$ for $M_{\eta,z}$ such that 
$$ \int {|p|}^2{|g|}^2 d\eta \leq \int {|p|}^2 d\theta$$ 
for every polynomial $p$, and $\eta|\partial\Omega_p$ is absolutely continuous with respect to $\theta$.\\

{\bf Lemma 4.9}.
Let $S$ be a cyclic subnormal tuple so that $S$ can be identified with $M_{\eta,z}$ for some compactly supported positive regular Borel measure $\eta$ on $\C^n$. Assume that supp$(\eta) \subset {\bar \Omega_p}$ and $\eta$ has no atoms on $\partial\Omega_p$. If there exists a bounded linear map $X: {P}^2(\eta) \rightarrow {P}^2(\theta)$ with dense range such that $XM_{\eta,z} = M_{\theta,z}X$, then there exists a cyclic vector $h$ for $M_{\theta,z}$ such that 
$$\int {|p|}^2{|h|}^2 d\theta \leq \int {|p|}^2 d(\eta|{\partial\Omega_p})$$
for every polynomial $p$, and $\theta$ is absolutely continuous with respect to $\eta|{\partial\Omega_p}$.\\

{\bf Proposition 4.10}.
Let $S$ be a cyclic subnormal tuple so that $S$ can be identified with $M_{\eta,z}$ for some compactly supported positive regular Borel measure $\eta$ on $\C^n$. Then $(S=) M_{\eta,z}$ is quasisimilar to $M_{\theta,z}$ if and only if\\
(a) there exists a cyclic vector $g$ for $M_{\eta,z}$ such that 
$$\int {|p|}^2{|g|}^2 d\eta \leq \int {|p|}^2 d\theta $$
for every polynomial $p$, and\\
(b) there exists a cyclic vector $h$ for $M_{\theta,z}$ such that 
$$\int {|p|}^2{|h|}^2 d\theta \leq \int {|p|}^2 d(\eta|{\partial\Omega_p})$$
for every polynomial $p$.\\

{\bf Proposition 4.11}.
Let $S$ be a cyclic subnormal tuple so that $S$ can be identified with $M_{\eta,z}$  for some compactly supported positive regular Borel measure $\eta$ on $\C^n$. Then $(S=) M_{\eta,z}$ is similar to $M_{\theta,z}$ if and only if there exist positive constants $c$ and $d$ such that
$$\int {|p|}^2 d\eta \leq c\int {|p|}^2 d\theta $$
and 
$$\int {|p|}^2 d\theta \leq d\int {|p|}^2 d(\eta|{\partial\Omega_p})$$
for every polynomial $p$. Also, $(S=) M_{\eta,z}$ is unitarily equivalent to $M_{\theta,z}$ if and only if $d\eta = {|h|}^2d\theta$ for some cyclic vector $h$ for $M_{\theta,z}$.\\

It would be interesting to know whether the statements of Propositions 4.10 and 4.11 remain valid even when $\theta$ has atoms on $\partial \Omega_p$. Since the surface area measure $\sigma_p$ on $\partial \Omega_p$ is not absolutely continuous with respect to the restriction $\nu_p|\partial\Omega_p$ of the volumetric measure $\nu_p$ to $\partial\Omega_p$, Lemma 4.9 shows in particular that $M_{\sigma_p,z}$  cannot be quasisimilar to $M_{\nu_p,z}$. This negative result can actually be extended to the multiplication tuples $M_{\sigma_p,z}$  and $M_{\nu_p,z}$ associated with the domains $\Sigma_p$. The next proposition generalizes \cite[Proposition 3.4 (d)]{At-Po} with an analogous proof; a complete proof is presented here for the reader's convenience. \\

{\bf Proposition 4.12}. There is no injective bounded linear map from $A^2(\nu_p)$ to $H^2(\sigma_p)$ that intertwines the  multiplication tuples $M_{\nu_p,z}$ and $M_{\sigma_p,z}$ associated with $\Sigma_p$.
\begin{proof}
We note that $\sum_{i=1}^n(M^*_{\sigma_p,z_i})^{p_{i,1}}(M_{\sigma_p,z_i})^{p_{i,1}}+\cdots+(M^*_{\sigma_p,z_i})^{p_{i,m_i}}(M_{\sigma_p,z_i})^{p_{i,m_i}}$ is the identity operator on $H^2(\sigma_p)$ so that $S\equiv ((M_{\sigma_p,z_1})^{p_{1,1}},\ldots,(M_{\sigma_p,z_n})^{p_{n,m_n}})$ is a spherical isometry. It follows from \cite[Proposition 2]{At2} that $S$ is subnormal and that the minimal normal extension $M$ of $S$ has its spectral measure $\rho_M$ supported on $\partial {\mathbb B}_Q$, where $Q =m_1+\cdots m_n$. It also follows from the Taylor functional calculus (see \cite{T}) and the spectral inclusion property for subnormal tuples (see \cite{Pu}) that the minimal normal extension $N$ of $T \equiv ((M_{\nu_p,z_1})^{p_{1,1}},\ldots,(M_{\nu_p,z_n})^{p_{n,m_n}})$ has its spectral measure $\rho_N$ supported on the closure $\bar {\mathbb B}_Q$ of ${\mathbb B}_Q$. Suppose there exists an injective  bounded linear map $Y: A^2(\nu_p)\rightarrow H^2(\sigma_p)$ satisfying $ YM_{\nu_p,z_i} =M_{\sigma_p,z_i}Y $ for all $i$. Then $Y$ also satisfies $YT_i = S_iY$ for all $i$. 
If $1_{\nu_p}$ is the constant function of  $A^2(\nu_p)$ taking value $1$, then for any $m$-variable polynomial $q \in \C[z]$ one has
$$ \int_{\partial {\mathbb B}_Q}|q(z)|^2d\|\rho_M(z)Y1_{\nu_p}\|^2 = \|q(S)Y1_{\nu_p}\|^2 = \|Yq(T)1_{\nu_p}\|^2 $$
$$\leq \|Y\|^2 \|q(T)1_{\nu_p}\|^2 = \|Y\|^2\int_{\bar {\mathbb B}_Q}|q(z)|^2d\|\rho_N(z)1_{\nu_p}\|^2.$$
Appealing to \cite[Theorem 3.5]{Ru}, we choose a sequence $\{q_n\}$ of polynomials in $\C[z]$ such that $q_n$ are bounded in absolute value by $1$, converge uniformly to $0$ on compact subsets of ${\mathbb B}_Q$, and satisfy
$$ \lim_{n \rightarrow \infty} |q_n(z)| = 1 \ \ z{\rm -a.e.\ } [\|\rho_M(\cdot)Y1_{\nu_p}\|^2].$$
Replacing $q$ by $q_n$ in the previous inequality, letting $n$ tend to infinity, and noting that the measure $\|\rho_N(\cdot)1_{\nu_p}\|^2$ vanishes on $\partial {\mathbb B}_Q$, we arrive at the absurdity $ 0 < \|Y1_{\nu_p}\|^2 \leq 0$.
\end{proof}

{\bf Remark 4.13}. Combining \cite[Theorem 2.3]{S} with our observation in the proof of Proposition 3.5 that the $\bar \partial$-Neumann operator $N_1$ correponding to $\Omega_p$ is compact, it follows that the short exact sequence of C$^*$-algebras 
$$ 0 \rightarrow {\mathcal K}(A^2(\nu_p)) \overset {\iota}\rightarrow {\A} \overset{\pi}\rightarrow C(\partial \Omega_p) \rightarrow 0$$
obtains, where ${\mathcal K}(A^2(\nu_p))$ and $\A$ are as in the proof of Proposition 3.7, $\iota$ is the inclusion map, and $\pi$ is a unital $*$-homomorphism satisfying $\pi(\tilde T_{\nu_p,\phi}) = \phi|{\partial \Omega_p}$ for any $\phi \in C(\bar \Omega_p)$. In view of Remark 4.2, even the $\bar \partial$-Neumann operator $N_1$ correponding to $\Sigma_p$ is compact; as such \cite[Theorem 2.3]{S} yields that the short exact sequence as recorded here obtains with $\Omega_p$ replaced by $\Sigma_p$ (and with the associated symbols interpreted accordingly). On the other hand, the short exact sequence of Remark 3.8 (b) was derived appealing to  \cite[Proposition 3.10]{Di-E-Ev} which necessitated that the multiplication tuple  $M_{\sigma_p,z}$ there be regular; this in turn forced us to use the full strength of the definition of $\Omega_p$ via Proposition 2.5. One may then ask in particular whether the short exact sequence of Remark 3.8 (b)  obtains with $\Omega_p$ replaced by $\Sigma_p$ - indeed, it does if $\Sigma_p$ is chosen to be a complex ellipsoid (see \cite[Theorem 2.1]{Cr-R}) and it also does if $\Sigma_p$ is chosen to be $\tilde \Sigma_p$ since a $\partial {\tilde \Sigma}_p$-isometry is an essntially normal $A(\tilde \Sigma_p)$-isometry (by Remark 4.2) and is moreover regular by the virtue of $\tilde \Sigma_p$ being strictly pseudoconvex (refer to \cite{E}).\\

While the main focus of the present paper has been on multivariable isometies associated with the domains $\Omega_p$, Proposition 3.3 as well as the analysis in the present section suggest that even subnormal tuples that have the spectral measures of their minimal normal extensions supported on $\bar \Omega_p$ (and not just on $\partial \Omega_p$) are worth exploring. To corroborate that assertion, we first proceed to verify that the domains $\Omega_p$ satisfy the properties (F1), (F2), (F3) and (F4) as enunciated in \cite[Section 1]{Di-E}. (It will also be clear that the domains $\Sigma_p$ satisfy the properties (F1), (F2) and (F4)).\\

(F1) The closure ${\overline {\Omega}_p}$ of $\Omega_p$ is a Stein compactum of $\C^n$: This follows from the fact that ${\overline {\Omega}_p}$ is a compact convex subset of $\C^n$ (refer to \cite[Chapter 3]{Ran}).\\

(F2) $\O({\overline {\Omega}_p})$, the vector space of functions $f$ such that $f$ is holomorphic on an open neighborhood $U_f$ of $\bar \Omega_p$, is weak*-dense in  $H^{\infty}(\Omega_p)$: This follows from our observation in the last line of Remark 4.3. \\

It may be recalled that $A(\Omega_p)$ is the closure of  $\O({\overline {\Omega}_p})$ in the sup norm with respect to $\bar \Omega_p$.\\

(F3) There exists a natural number $N$ and and an injective mapping $f \in A(\Omega_p)^N$ such that the image of the Shilov boundary of $A(\Omega_p)$ is contained in the topological boundary of the unit ball ${\mathbb B}_{N}$: This follows from Proposition 2.4 and our observations in the proof of Proposition 2.5. \\

(F4) There exists a positive regular Borel measure $\mu$ supported on the Shilov boundary of $\Omega_p$ (which, as we know, is $\partial \Omega_p$) such that the canonical map $r_{\mu}$ from $\O({\overline {\Omega}_p}) \rightarrow L^{\infty}(\mu)$ extends to an algebra homomorphism $\tilde r_{\mu} : H^{\infty}(\Omega_p) \rightarrow L^{\infty}(\mu)$ that is isometric and weak$^*$-continuous (which is the same as calling $\mu$ a `faithful Henkin measure'): Since the non-tangential boundary limit of any $f \in \O({\overline {\Omega}_p})$ is the restiction of $f$ to $\partial \Omega_p$, the normalized surface area measure $\sigma_p$ on $\partial \Omega_p$ is a faithful Henkin measure in the light of Remark 4.3.\\

{\bf Remark 4.14}. The preceding observations allow us to apply \cite[Theorem 1.4]{Di-E} to those operator tuples $T \in {\mathcal B}({\mathcal H})^n$ that possess an isometric and a weak$^*$-continuous $H^{\infty}(\Omega_p)$-functional calculus $\Phi_T: H^{\infty}(\Omega_p) \rightarrow {\mathcal B}({\mathcal H})$ (satisfying $\Phi_T(1) = I$ and $\Phi_T(z_i) = T_i$ for all $i$) so that, for such tuples $T$, we have the following: The weak operator topology and the weak$^*$ operator topology coincide on the algebra $\Phi_T(H^{\infty}(\Omega_p))$ and any unital weak$^*$-closed subalgebra of $\Phi_T(H^{\infty}(\Omega_p))$ is reflexive (cf. Remark 2.7 (a)).\\

let $T \in {\mathcal B}({\mathcal H})^n$ be an operator tuple possessing a contractive and a weak$^*$-continuous $H^{\infty}(\Omega_p)$-functional calculus $\Phi_T$. Suppose further that $T$ has its Taylor spectrum {\it dominating in $\Omega_p$} so that the sup norm of any $f \in H^{\infty}(\Omega_p)$ equals the supremum of $|f|$ over the intersection of $\Omega_p$ with the Taylor spectrum $\sigma(T)$ of $T$.  Since $\Omega_p$ is a bounded convex domain with smooth boundary, $\Omega_p$ satisfies the `Gleason property' so that, for any $a \in \Omega_p$ and any $f \in H^{\infty}(\Omega_p)$, one has
$$ f(z) -f(a) = \sum_{i=1}^n (z_i-a_i)f_i(z) \ \ \ (z \in \Omega_p),$$
where the so-called Leibenzon divisors $f_i$ are given by 
$$f_i(z) = \int_0^1\frac{\partial f}{\partial z_i}(a+t(z-a))dt $$
and are in $H^{\infty}(\Omega_p)$ (refer to \cite{Gr} and \cite{He}).
Using this and arguing exactly as in \cite[Lemma 2.3.6]{Di1} one can prove that, for any $f \in H^{\infty}(\Omega_p)$, $f(\sigma(T) \cap \Omega_p)$ is contained in the Taylor spectrum of $\Phi_T(f)$; that easily leads to the sup norm of $f$ with respect to $\Omega_p$ being less than or equal to $\|\Phi_T(f)\|$. Thus, in this case,  the functional calculus $\Phi_T$ is indeed isometric.\\


\newpage

\begin{thebibliography}{00}

\bibitem{A} A.B. Aleksandrov,  Inner functions on compact spaces, Funct. Anal. Appl. 18 (1984) 87-98.

\bibitem{At1} A. Athavale, Subnormal tuples quasi-similar to the Szeg\"o tuple, Michigan Math. J. 35 (1988) 409-412.

\bibitem{At2} A. Athavale,  On the intertwining of joint isometries, J. Operator Theory 23 (1990) 339-350.

\bibitem{At3} A. Athavale,  On the intertwining of $\partial {\mathcal D}$-isometries, Complex Anal. Oper. Theory 2 (2008) 417-428. 

\bibitem{At-P} A. Athavale and S. Pedersen,  Moment problems and subnormality, J. Math anal. Appl. 146 (1990) 434-441.

\bibitem{At-Po} A. Athavale and S. Podder, on the multiplication tuples related to certain reproducing kernel Hilbert spaces, Complex Anal. Oper. Theory 10 (2016) 1329-1338.

\bibitem{Bo-S} H.P. Boas and E.J. Straube,  Global regularity of the $\bar \partial$-Neumann problem: A survey of the $L^2$-Sobolev theory, Several Complex Variables, MSRI Publications 37 (1999) 79-111.

\bibitem{Car} E. Cartan, Sur les domaines born\'es homog\`enes de l'espace des $n$ variables complexes, Abh. Math. Sem. Hamburg 11 (1935) 116-162.

\bibitem{Cat} D.W. Catlin,  Global regularity of the $\bar \partial$-Neumann problem, Proceedings of Symposia in Pure  Mathematics, Vol. 4 (1984) 39-49.

\bibitem{Co1} J.B. Conway,  The theory of subnormal operators, Math. Surveys and Monographs, Vol. 36, AMS, Providence, R.I., 1991.

\bibitem{Co2} J.B. Conway, A course in operator theory, Graduate Studies in Mathematics, Vol. 21, AMS, Providence, R.I., 2000.

\bibitem{Cr-R} D. Crocker and I. Raeburn, Toeplitz operators on certain weakly pseudoconvex domains, J. Austral. Math. Soc. (Series A) 31 (1981) 1-14.

\bibitem{Cu} R.E. Curto, Spectral inclusion for doubly commuting subnormal $n$-tuples, Proc. Amer. Math. Soc. 83 (1981) 730-734.

\bibitem{D} K.R. Davidson,  On operators commuting with Toeplitz operators modulo the compact operators, J. Funct. Anal. 24 (1977) 356-368.

\bibitem{Di1} M. Didas, Dual algebras generated by von Neumann $n$-tuples over strictly pseudoconvex sets, Dissertationes Math. (Rozprawy Mat.) 425 (2004).

\bibitem{Di2} M. Didas,  A note on the Toeplitz projection associated with spherical isometries, preprint.

\bibitem{Di-E} M. Didas and J. Eschmeier, Subnormal tuples on strictly pseudoconvex and bounded symmetric domains, Acta Sci. Math. (Szeged) 71 (2005) 691-731.

\bibitem{Di-E-Ev} M. Didas, J. Eschmeier and K. Everard,  On the essential commutant of analytic Toeplitz operators associated with spherical isometries, J. Funct. Anal. 261 (2011) 1361-1383.

\bibitem{Die-F} K. Diedrich and J.E. Fornaess,  Pseudoconvex domains with real analytic boundary, Ann. Math. (1978) 371-384.

\bibitem{E} J. Eschmeier,  On the reflexivity of multivariable isometries, Proc. Amer. Math. Soc. 134 (2005) 1783-1789.

\bibitem{E-Ev} J. Eschmeier and K. Everard, Toeplitz projections and essential commutants, J. Funct. Anal. 269 (2015) 1115-1135.

\bibitem{F-Ko} G.B. Folland and J.J. Kohn,  The Neumann problem for the Cauchy-Riemann complex, Annals of Mathematics Studies, Princeton University Press, Princeton, New Jersey, 1972.


\bibitem{Gr} M. Grang\'{e},  Diviseurs de Leibenson et probl\`{e}me de Gleason pour $H^{\infty}(\Omega)$ dans le cas convexe, Bull. Soc. Math. France 114 (1986) 225-245.

\bibitem{H-S} M. Hakim and N. Sibony, Fronti\'ere de Shilov et spectre de $A(\bar D)$ pour les domaines faiblement psedoconvexes, C. R. Acad. Sci. Paris 281 (1975) 959-962.

\bibitem{Ha} W.W. Hastings,  Commuting subnormal operators simultaneously quasisimilar to unilateral shifts, Illinois J. Math. 22 (1978) 506-519.

\bibitem{He} G. M. Henkin,  The approximation of functions in pseudo-convex domains and a theorem of Z. L. Leibenzon (Russian), Bull. Acad. Polon. Sci. S\'{e}r. Sci. Math. Astronom. Phys. 19 (1971) 37-42.

\bibitem{Ho} L. H\"ormander, $L^2$ estimates and existence theorems for the $\bar \partial$ operator, Acta Math. 113 (1965) 89-152.

\bibitem{I} T. Ito, On the commuting family of subnormal operators, J. Fac. Sci. Hokkaido Univ. 14 (1958) 1-15.

\bibitem{J-P} M. Jarnicki and P. Pflug,  First steps in several complex variables: Reinhardt domains, European Mathematical Society, 2008.

\bibitem{Je-L} N.P. Jewell and A.R. Lubin,  Commuting weighted shifts and analytic function theory in several variables, J. Operator Theory 1 (1979) 207-223.

\bibitem{K-Pt} K. Kli\'s and M. Ptak,  $k$-hyperreflexive subspaces, Houston J. Math. 32 (2206) 299-313.

\bibitem{Ko} J.J. Kohn, The range of the tangential Cauchy-Riemann operator, Duke Math. J. 53 (1986) 525-545.  

\bibitem{Ko-R} J.J. Kohn and H. Rossi,  On the extension of holomorphic functions from the boundary of a complex manifold, Ann. Math. 81 (1965) 451-472.

\bibitem{Kr} S. Krantz,  Function theory of several complex variables, AMS, Providence, R.I., 2001.


\bibitem{M}W. Mlak, Intertwining operators, Studia Math. 43 (1972) 219-233.


\bibitem{Pf} P. Pflug, \"Uber polynomiale funktionen auf holomorphiegebieten, Math. Z. 138 (1974) 133-139.

\bibitem{Pi} S.I. Pinchuk, Homogeneous domains with piecewise-smooth boundaries, Mat. Zametki 32 (1982); English transl.  Math Notes 32 (1982) 849-852.

\bibitem{Pr} B. Prunaru, Some exact sequences for Toeplitz algebras of spherical isometries, Proc. Amer. Math. Soc. 135 (2007) 3621-3630.

\bibitem{Pu} M. Putinar, Spectral inclusion for subnormal $n$-tuples, Proc. Amer. Math. Soc. (1984) 405-406.

\bibitem{Ra-St} Raich and Straube, Compactness of the complex Green operator, Math. Res. Lett. 15 (2008) 761-778.

\bibitem{Ran} R. M. Range,  Holomorphic functions and integral representations in several complex variables, Springer-Verlag, New York, 1986.

\bibitem{Ru} W. Rudin, New constructions of functions holomorphic in the unit ball of $\C^n$, CBMS Regional Conference Series in Mathematics, AMS, Providence, R.I., 1986.

\bibitem{S} N. Salinas, The $\bar \partial$-formalism and the C$^*$-algebra of the Bergman $n$-tuple, J. Operator Theory 22 (1989) 325-343.

\bibitem{S-Sh-U} N. Salinas, A. Sheu and H. Upmeier, Toeplitz Operators on Pseudoconvex Domains and Foliation C*-algebras, Ann. Math. 130 (1989) 531-565.

\bibitem{Sh} A. Sheu, Isomprphism of the Toeplitz C$^*$-algebras for the Hardy and Bergman spaces of certain Reinhardt domains, Proc. Amer. Math. Soc. 116 (1992) 113-120.

\bibitem{Su} T. Sunada, Holomorphic equivalence problem for bounded Reinhardt domains, Math. Ann. 235 (1978) 111-128.

\bibitem{T} J. L. Taylor, The analytic-functional calculus for several commuting operators, Acta Math. 125 (1970) 1-38.

\bibitem{U} H. Upmeier, Toeplitz operators and index theory in several complex variables, Operator Theory Advances and Applications, Vol. 81, Birkh\"auser Verlag, Basel, 1996.

\end{thebibliography}
\end{document}